\documentclass[12pt,leqno]{article}
\usepackage{amsmath, amssymb}

\def\epsilon{\varepsilon}

\begin{document}

\LARGE \noindent 
{\bf A remark on locally direct product subsets
in a topological Cartesian space}

\large

\vspace*{0.4em}

\hfill Hiroki Yagisita (Kyoto Sangyo University)

\vspace*{1.2em}

\normalsize

\noindent 
Abstract: 

Let $X$ and $Y$ be topological spaces. 
Let $C$ be a path-connected closed set of $X\times Y$. 
Suppose that $C$ is locally direct product, that is, 
for any $(a,b)\in X\times Y$, there exist 
an open set $U$ of $X$, an open set $V$ of $Y$, 
a subset $I$ of $U$ and a subset $J$ of $V$ 
such that $(a,b) \in U\times V$ and 
$$C\cap (U\times V)=I\times J$$
hold. Then, in this memo, we show that $C$ is globally so, 
that is, there exist 
a subset $A$ of $X$ and a subset $B$ of $Y$ 
such that 
$$C=A\times B$$
holds. The proof is elementary. 
Here, we note that one might be able to think 
of a (perhaps, open) similar problem 
for a fiber product of locally trivial fiber spaces, 
not just for a direct product of topological spaces. 

In Appendix, we mentioned a simple example 
of a $C([0,1];\mathbb R)$-manifold that cannot be embedded 
in the direct product $(C([0,1];\mathbb R))^n$ 
as a $C([0,1];\mathbb R)$-submanifold. 
In addition, we introduce the concept of topological 2-space, 
which is locally the direct product of topological spaces
and an analog of homotopy category for topological 2-space. 
Finally, we raise a question on the existence of an $\mathbb R^n$-Morse function 
and the existence of an $\mathbb R^n$-immersion 
in a finite-dimensional $\mathbb R^n$-Euclidean space. 
Here, we note that the problem of defining the concept 
of an $\mathbb R^n$-handle body may also be considered. 

\vspace*{0.8em} 
\noindent 
Keywords: \ 
Morse-Bott index,  compact foliation, fiber bundle, fibration,
holomorphic web, Kodaira-Spencer map, Kuranishi family, Stein manifold.  

\vfill 
\vspace*{0.8em} 
\noindent 
The related literature: \ 
``https://www.researchgate.net/profile/Hiroki\_Yagisita''

\noindent

\noindent 


\noindent 

\newpage
\noindent 
{\bf Proof} : \ 

\noindent 
[Step 1] \ \ \ 

Let $(x,y)$ is a continuous mapping from $[0,1]$ to $C$. 
Then, in this step, we show that $(x(0),y(1)) \in C$ and $(x(1),y(0)) \in C$ hold. 

Let $$T \, := \, \{ \, t\in [0,1] \, | \, 
\forall \, (u,v) \in [0,t]\times[0,t]: \, 
(x(u),y(v)) \in C \, \}.$$
Then, $0\in T$ holds. 
Because $C$ is a closed set of $X\times Y$, 
$T$ is a closed set of $[0,1]$.

We show that for any $t_0\in [0,1)\cap T$, there exists $t_1 \in (t_0,1]$ 
such that 
$[0,t_1] \subset T$ holds. 
Let 
$$t_0\in [0,1)\cap T.$$
Let 
$$D \, := \, \{ \, (u,v) \in [0,1]\times [0,1] \, | \, 
(x(u),y(v)) \in C \, \}.$$
Then, $D$ is a closed set of $[0,1]\times [0,1]$ and 
$$\cup_{t\in[0,1]}\{(t,t)\} \subset D$$
and 
$$[0,t_0]\times[0,t_0]\subset D$$ hold. 
Because $C$ is a locally direct product set of $X\times Y$, 
$D$ is a locally direct product set of $[0,1]\times [0,1]$. 
That is,  for any $(a,b)\in [0,1]\times [0,1]$, there exist 
open sets $U$ and $V$ of $[0,1]$, 
a subset $I$ of $U$ and a subset $J$ of $V$ 
such that $(a,b) \in U\times V$ and 
$$D\cap (U\times V)=I\times J$$
hold. So, there exist open sets $U_0$ and $V_0$ of $[0,1]$, 
a subset $I_0$ of $U_0$ and a subset $J_0$ of $V_0$ 
such that $(t_0,t_0) \in U_0\times V_0$ and 
$$D\cap (U_0\times V_0)=I_0\times J_0$$
hold. Then, from $t_0\in [0,1)$, there exists $t_2\in (t_0,1]$ such that 
$$[t_0,t_2]\times [t_0,t_2] \subset U_0 \times V_0$$
holds. So, from $\cup_{t\in[t_0,t_2]}\{(t,t)\} \subset D$, 
$\cup_{t\in[t_0,t_2]}\{(t,t)\} \subset I_0\times J_0$ holds. 
Hence, $[t_0,t_2]\times [t_0,t_2] \subset I_0 \times J_0$ holds. So,  
\newpage
$$[t_0,t_2]\times [t_0,t_2] \subset D$$
holds. Well, there exist $n\in \mathbb N$, $a_1, a_2, \cdots, a_n \in [0,t_0]$, 
open sets $U_1, V_1, U_2, V_2, \cdots, U_n, V_n$ of $[0,1]$ 
and sets $I_1, J_1, I_2, J_2, \cdots, I_n, J_n$ such that 
$$[0,t_0]\times\{t_0\} \subset \cup_{k=1,2,\cdots,n} (U_k\times V_k)$$
holds and for any $k\in\{1,2,\cdots,n\}$, $I_k\subset U_k$, $J_k\subset V_k$, 
$(a_k,t_0)\in U_k\times V_k$ and 
$$D\cap (U_k\times V_k)=I_k\times J_k$$
hold. Then, from $0\leq t_0<t_2\leq1$, there exists $t_3\in (t_0,t_2]$ 
such that
$$[t_0,t_3] \subset \cap_{k=1,2,\cdots,n} V_k$$ 
holds. Let 
$$S \, := \, \{ \, u\in[0,t_0] \, | \, [u,t_0] \times [t_0,t_3] \subset D\, \}.$$
Then, from $0\leq t_0<t_3\leq t_2\leq 1$, $t_0\in S$ holds. 
Because $D$ is a closed set of $[0,1]\times [0,1]$, 
$S$ is a closed set of $[0,t_0]$.  
Now, we show that for any $u_0 \in (0,t_0]\cap S$, there exists 
$u_1 \in [0,u_0)$ such that 
$[u_1,t_0] \subset S$ holds. 
Let $u_0 \in (0,t_0]\cap S$. 
Then, there exists $k_0\in\{1,2,\cdots,n\}$ such that 
$(u_0,t_0)\in U_{k_0}\times V_{k_0}$ holds. 
From $u_0 \in (0,t_0]$, there exists 
$u_1 \in [0,u_0)$ such that 
$[u_1,u_0] \subset U_{k_0}$ holds. 
Because $t_0\in V_{k_0}$ and $[u_1,u_0]\times\{t_0\}\subset D$ hold, 
$[u_1,u_0]\times\{t_0\}\subset I_{k_0}\times J_{k_0}$ holds. So, 
$$[u_1,u_0]\subset I_{k_0}$$
holds. On the other hand, from $u_0\in S$ and $[t_0,t_3]\subset V_{k_0}$, 
$\{u_0\}\times [t_0,t_3]\subset D\cap (U_{k_0}\times V_{k_0})=I_{k_0}\times J_{k_0}$ 
holds. So, 
$$[t_0,t_3]\subset J_{k_0}$$
holds. Therefore, $[u_1,u_0]\times[t_0,t_3]\subset D$ holds. 
Hence, $u_1\in S$ holds. So, for any $u_0 \in (0,t_0]\cap S$, there exists 
$u_1 \in [0,u_0)$ such that 
$[u_1,t_0] \subset S$ holds. 
Hence, $S$ is a nonempty closed open set of $[0,t_0]$. 
Because $S=[0,t_0]$ holds, 
$$[0,t_0] \times [t_0,t_3] \subset D$$
holds. Similarly, there exists $t_4\in (t_0,t_2]$ such that 
\newpage
$$[t_0,t_4] \times [0,t_0] \subset D$$
holds. Therefore, as we set $t_1:=\min\{t_3,t_4\}$, 
$t_1\in (t_0,1]$ and $[0,t_1] \times [0,t_1] \subset D$ hold. 
So, for any $t_0\in [0,1)\cap T$, there exists $t_1 \in (t_0,1]$ 
such that 
$[0,t_1] \subset T$ holds. 

Therefore, $T$ is a nonempty closed open set of $[0,1]$. 
Because $T=[0,1]$ holds, $(x(0),y(1)) \in C$ and $(x(1),y(0)) \in C$ hold. 
\hfill ---

\noindent 
[Step 2] \ \ \ 

Let 
$$A \, := \, \{ \, x\in X \, | \, \exists \, y\in Y: (x,y) \in C \, \}$$
and 
$$B \, := \, \{ \, y\in Y \, | \, \exists \, x\in X: (x,y) \in C \, \}.$$
Then, $C\subset A\times B$ holds. Let $(x_0,y_1)\in A\times B$. 
Then, we show that $(x_0,y_1)\in C$ holds. There exist 
$y_0\in Y$ and $x_1\in X$ such that 
$(x_0,y_0)\in C$ and $(x_1,y_1)\in C$ hold. 
Because $C$ is path-connected, in virtue of Step 1, 
$(x_0,y_1)\in C$ holds. 
\hfill 
$\blacksquare$

\vspace*{0.8em}



\noindent 
{\bf Comment} : \ 
Does there exist a $\mathbb C^2$-manifold $N$ 
such that for any $\mathbb C$-manifolds 
$M_1$ and $M_2$, $N$ can not be embedded 
in $M_1\times M_2$ as a $\mathbb C^2$-submanifold ? 
Kasuya proposed a candidate for such a compact $\mathbb C^2$-manifold $N$.
Our result may be useful to prove that it is such. 
For the definition of a $\mathbb C^n$-manifold, see [1] or [2]. 

In [1], we proposed a candidate of a noncompact $\mathbb C^2$-manifold $N$ 
such that (1) for any $\mathbb C$-manifolds 
$M_1$ and $M_2$, $N$ can not be embedded 
in $M_1\times M_2$ as a $\mathbb C^2$-submanifold 
but (2) there exists $k\in \mathbb N$ such that 
$N$ can be embedded 
in the $k$-dimensional Euclidean $\mathbb R^2$-space $(\mathbb R^k)^2$ 
as an $\mathbb R^2$-submanifold.  
\hfill ---


\noindent 
{\bf Remark} : \ 
For deformation of $\mathbb R^n$-structures 
and $\mathbb C^n$-structures, 
see Kodaira and Spencer, 
{\it Ann. of Math.}, 74 (1961), 52-100. 
\hfill ---

\noindent
{\bf Problem} : \ 
A $\mathbb C$-holomorphically convex domain of $\mathbb C^n$ 
is an $n$-dimensional Stein manifold. 
So, it is a closed $\mathbb C$-submanifold of $\mathbb C^{2n+1}$. 
Let $C$ be a $\mathbb C^2$-holomorphically convex domain. 
Then, is $C$ the direct product of some Stein manifolds ? 
For the definition of $\mathbb C^n$-holomorphic convexity, see [2]. 
\hfill ---




\vfill
\vspace*{0.8em}

\noindent  
Acknowledgment: \ \ \ 
As in Comment, Professor Naohiko Kasuya proposed it. 
This work was supported by JSPS KAKENHI Grant Number JP16K05245. 

\newpage

\noindent 
{\large \bf Appendix 0}

\vspace*{0.4em}

In [3], we gave a simple example 
of a connected metrizable $1$-dimensional
$C([0,1];\mathbb R)$-manifold
that cannot be embedded 
in the direct product space $(C([0,1];\mathbb R))^n$ 
as a $C([0,1];\mathbb R)$-submanifold. 

\vspace*{1.6em}

\noindent 
{\large \bf Appendix 1}

\vspace*{0.4em}

\noindent 
{\bf Definition 1} (2-space) : 

$W:=(W,S)$ is said to be a 2-space and 
$S$ is said to be the system of local Cartesian neighborhoods
of $W$, if it satisfies the followings. 

(1) \ $W$ is a topological space. $S$ is a set. 

(2) \ Let $\varphi \in S$. Then, $\varphi$ is a homeomorphism 
from an open set $W_\varphi$ of $W$ to the direct product 
of topological spaces $U_\varphi$ and $V_\varphi$.

(3) \ Let $\varphi_1, \varphi_2 \in S$ and $c\in W_{\varphi_1}\cap W_{\varphi_2}$. Then, 
there exist an open set $M_1$ of $U_{\varphi_1}$, an open set $N_1$ of $V_{\varphi_1}$, 
a map $f$ from $M_1$ to $U_{\varphi_2}$ and a map $g$ from $N_1$ to $V_{\varphi_2}$ 
such that $c\in \varphi_1^{-1}(M_1\times N_1) \subset  W_{\varphi_2}$ holds 
and for any $w\in \varphi_1^{-1}(M_1\times N_1) $, 
$$\varphi_2(w)=(f(\pi_U(\varphi_1(w))),  g(\pi_V(\varphi_1(w))) )$$
holds. Here,  $\pi_U$ is the projection from $U_{\varphi_1}\times V_{\varphi_1}$ 
to $U_{\varphi_1}$ and $\pi_V$ is the projection from $U_{\varphi_1}\times V_{\varphi_1}$ 
to $V_{\varphi_1}$.

(4) \ $W \, = \, \cup_{\, \varphi\in S} \, W_\varphi$ holds.

\noindent 
{\bf Example 2} (locally direct product subset) : 

Let $X$ and $Y$ be topological spaces. Then, 
a locally direct product subset of $X\times Y$ is a 2-space. 

\noindent 
{\bf Example 3} (2-product) : 

Let $(W_1,S_1)$ and $(W_2,S_2)$ be 2-spaces. 
Then, for any $(\varphi_1,\varphi_2)\in S_1\times S_2$, 
the map $$(\varphi_1, \varphi_2) \, : \, W_{1,\varphi_1}\times W_{2,\varphi_2}
\longrightarrow (U_{1,\varphi_1}\times U_{2,\varphi_2})
\times (V_{1,\varphi_1}\times V_{2,\varphi_2})$$
is a homeomorphism. So, the 2-product 
$$W_1\times_2 W_2 \, := \, (W_1\times W_2, S_1\times S_2)$$
is a 2-space. 
\hfill ---

\newpage 

\noindent 
{\bf Definition 4} (2-map) : 

Let $h$ be a continuous map from a 2-space $(W_1,S_1)$ to a 2-space $(W_2,S_2)$. 
Then, $h$ is said to be a 2-map, if it satisfies the following. 

Let $c\in W_1$, $\varphi_1\in S_1$, $c\in W_{1,\varphi_1}$, 
$\varphi_2\in S_2$ and $h(c)\in W_{2,\varphi_2}$. Then,  
there exist an open set $M_1$ of $U_{1,\varphi_1}$, an open set $N_1$ of $V_{1,\varphi_1}$, 
a map $f$ from $M_1$ to $U_{2,\varphi_2}$ and a map $g$ from $N_1$ to $V_{2,\varphi_2}$ 
such that $c\in \varphi_1^{-1}(M_1\times N_1)\subset  h^{-1}(W_{2,\varphi_2})$ holds  
and for any $w\in \varphi_1^{-1}(M_1\times N_1) $, 
$$\varphi_2(h(w))=(f(\pi_U(\varphi_1(w))),  g(\pi_V(\varphi_1(w))) )$$
holds. Here,  $\pi_U$ is the projection from $U_{1,\varphi_1}\times V_{1,\varphi_1}$ 
to $U_{1,\varphi_1}$ and $\pi_V$ is the projection from $U_{1,\varphi_1}\times V_{1,\varphi_1}$ 
to $V_{1,\varphi_1}$. 

\noindent 
{\bf Definition 5} (2-homotopy) : 

Let $h_0$ and $h_1$ be 2-maps from a 2-space $W_1$ 
to a 2-space $W_2$. 
Then, $h_0$ and $h_1$ said to be 2-homotopic, 
if there exists a continuous map $H$ from $[0,1]\times W_1$ to $W_2$ 
such that for any $t\in [0,1]$, the map $w_1 \in W_1 \, \mapsto \, H(t,w_1) \in W_2$ 
is a 2-map and for any $w_1 \in W_1$, $H(0,w_1)=h_0(w_1)$ 
and $H(1,w_1)=h_1(w_1)$ hold.  

\noindent 
{\bf Example 6} : 

2-spaces $\{0\}\times (\mathbb R/\mathbb Z)$ 
and $(\mathbb R/\mathbb Z)\times \{0\}$ 
are homeomorphic. However, they are not 2-homotopy equivalent.
\hfill ---

\noindent 
{\bf Problem} : 

(1) \ For a compact Hausdorff space $T$, 
define $T$-space, $T$-product, $T$-map and $T$-homotopy. 

(2) \ Do there exist a contractible compact Hausdorff space $T$ 
and $T$-spaces such that the $T$-spaces are homeomorphic 
but not $T$-homotopy equivalent ? 
\mbox{} \hfill ---

\noindent 
{\large \bf Appendix 2}

\vspace*{0.4em}

Let $M$ be a paracompact connected $\mathbb R^n$-manifold. 
Let $(f_1,f_2,\cdots,f_n)$ be an $\mathbb R^n$-map 
from $M$ to $\mathbb R^n$. 
Suppose that for any local Cartesian neighborhood 
$U_1\times U_2\times \cdots \times U_n$,  
$(a_1,a_2,\cdots,a_n)\in U_1\times U_2\times \cdots \times U_n$ 
and $k\in \{1,2,\cdots,n\}$, the map 
$$x_k \, \in \, U_k \ \ \ \mapsto \ \ \ 
f_k(a_1,a_2,\cdots,a_{k-1},x_k,a_{k+1},a_{k+2},\cdots,a_n) \, \in \, \mathbb R$$
is a Morse function on $U_k$. 
Then, does there exist $m\in \mathbb N$ 
such that $M$ can be $\mathbb R^n$-immersed  
in the $m$-dimensional $\mathbb R^n$-affine space 
$(\mathbb R^m)^n$ ?

\noindent 
{\bf Related problem} :  \ Let $M$ be a compact foliation. Then, does there exist 
a Morse-Bott function $f$ such that 
a critical manifold of $f$ is a leaf of $M$ 
and a leaf of $M$ 
is contained in a level set of $f$ ?
\hfill ---

\newpage 

\noindent 
{\bf References}

[1] H. Yagisita, Holomorphic differential forms of complex manifolds on commutative Banach algebras and a few related problems, {\it preprint}.

[2] H. Yagisita, Cartan-Thullen theorem for a $\mathbb C^n$-holomorphic convexity 
and a related problem, {\it preprint}.

[3] H. Yagisita, A manifold on the real commutative Banach algebra $C([0,1];\mathbb R)$
that cannot be embedded in the finite-dimensional Euclidean space
$C([0,1];\mathbb R^n)$, {\it preprint}. 

[4] H. Yagisita, Finite-dimensional complex manifolds on commutative Banach algebras and continuous families of compact complex manifolds, 
{\it Complex Manifolds}, 6 (2019), 228-264. 

\vfill 
\vspace*{0.8em} 
\noindent 
The related literature: \ 
``https://www.researchgate.net/profile/Hiroki\_Yagisita''


\end{document}